\documentclass[letterpaper, 10 pt, conference]{ieeeconf}  

\IEEEoverridecommandlockouts                   
\usepackage{cite}
\usepackage{amsmath,amssymb,amsfonts}

\usepackage{graphicx}
\graphicspath{ {./images/} }
\usepackage{comment}
\usepackage{enumerate}
\usepackage{mathtools}
\usepackage{soul}
\usepackage{color}
\usepackage{tikz}
\usetikzlibrary{arrows.meta,positioning,fit,calc}

\usepackage{etoolbox}

\BeforeBeginEnvironment{theorem}{\vspace{2mm}}
\AfterEndEnvironment{theorem}{\vspace{2mm}}

\BeforeBeginEnvironment{lemma}{\vspace{2mm}}
\AfterEndEnvironment{lemma}{\vspace{2mm}}

\BeforeBeginEnvironment{remark}{\vspace{2mm}}
\AfterEndEnvironment{remark}{\vspace{2mm}}

\BeforeBeginEnvironment{proposition}{\vspace{2mm}}
\AfterEndEnvironment{proposition}{\vspace{2mm}}

\BeforeBeginEnvironment{definition}{\vspace{2mm}}
\AfterEndEnvironment{definition}{\vspace{2mm}}

\BeforeBeginEnvironment{assumption}{\vspace{2mm}}
\AfterEndEnvironment{assumption}{\vspace{2mm}}

\BeforeBeginEnvironment{corollary}{\vspace{2mm}}
\AfterEndEnvironment{corollary}{\vspace{2mm}}

\title{\LARGE \bf
MPC and System Identification with Differentiable Physics: Fluid System and Particle Beam Control 
}

\author{Alan  Williams$^{1}$ \quad  Alp M. Sunol$^{2}$
\thanks{$^{1}$Alan Williams is with the Accelerator Operations and Technology - Instrumentation and Control (AOT-IC) group, Adaptive Machine Learning Team at Los Alamos National Lab, Los Alamos, NM 87545, USA (e-mail: awilliams@lanl.gov).}
\thanks{$^{2}$Alp M. Sunol is in the Brenner Group at Harvard University's School of Engineering and Applied Sciences (e-mail: asunol@seas.harvard.edu).}
}

\begin{document}

\maketitle
\thispagestyle{empty}
\pagestyle{empty}

%%%%%%%%%%%%%%%%%%%%%%%%%%%%%%%%%%%%%%%%%%%%%%%%%%%%%%%%%%%%%%%%%%%%%%%%%%%%%%%%
\begin{abstract}
We consider the problem of simultaneous control and parameter estimation when the model is available only as a differentiable physics simulator. We propose a receding-horizon control framework in which a model predictive control (MPC) objective is optimized using gradients obtained by differentiating through the simulator, while physical parameters are updated online using measurement data. Unlike classical MPC, which relies on explicit algebraic models, our approach treats the dynamics as a computational object and performs simulation-based optimization using automatic differentiation. A shared differentiable model enables joint, real-time optimization of control inputs and physical parameters. We present two preliminary examples to demonstrate the proposed framework on two challenging applications: a fluid flow problem and a particle accelerator.
\end{abstract}

%%%%%%%%%%%%%%%%%%%%%%%%%%%%%%%%%%%%%%%%%%%%%%%%%%%%%%%%%%%%%%%%%%%%%%%%%%%%%%%%
\section{Introduction}

Model-based control methods rely fundamentally on the availability of a system model that can be used to synthesize control laws. In classical control settings, this model is typically expressed in explicit algebraic form, for instance, as a linear or nonlinear state-space representation, which enables analytical reasoning, linearization, and the formulation of tractable optimization problems. In many modern applications, however, high-fidelity system models are no longer available in closed form. Instead, the system dynamics are encoded in complex computational pipelines, such as physics-based simulators, numerical solvers, or multiphysics models, which can only be queried through forward evaluation. This raises a fundamental question: how should control be performed when the system model is available only as a computational object rather than an explicit set of equations?

Recent advances in differentiable physics and automatic differentiation have made it possible to compute exact gradients of simulator outputs with respect to inputs, initial conditions, and physical parameters. As a result, such computational models can now be treated as differentiable maps, enabling gradient-based optimization through time-evolving simulations. This capability opens the door to a new class of model-based control methods, in which control inputs are optimized directly through simulation rollouts, rather than through explicit algebraic reformulations of the dynamics. At the same time, the availability of simulator gradients naturally supports online parameter estimation, allowing physical parameters within the model to be adapted using measurement data.

In this work, we formalize this approach as a receding-horizon framework that couples control and estimation through a shared differentiable model, treating the dynamics as an implicit computational constraint enforced through simulation.

%but also introduces challenges related to identifiability and persistence of excitation.

\section{Background and Literature} \label{sec:background_and_lit}

%\section{Background and Literature} \label{sec:background_and_lit}
\textit{Differentiable simulators} provide gradients of simulation outputs with respect to inputs, states, and physical parameters, enabling gradient-based optimization through complex computational models; see \cite{newbury2024review} for a recent review. When MPC is used with physics simulators, repeated rollouts can make gradient computation expensive, and methods such as \cite{liang2025robust} seek to reduce this cost. Only a limited number of works consider MPC with differentiable physics. In \cite{chen2022real}, the authors couple MPC (solved with iLQR) and online system identification using a differentiable rigid-body simulator, together with a heuristic certainty score derived from closed-form dynamics to detect insufficient excitation. In contrast, we assume the dynamics are available only through a simulator and focus on formulating MPC, state estimation, and online parameter estimation for higher-complexity GPU-based simulations.

This work presents two examples of differentiable simulation for control and parameter estimation: particle accelerator systems and fluid systems. In accelerator applications, model-based control has the potential to provide substantial speedups in tuning, since current approaches are often either model-free~\cite{scheinker2013model, williams2024experimental} or based on data-driven surrogate models~\cite{edelen2016neural, hirlaender2020model}. Neither approach makes efficient use of high-fidelity computational models, and both typically require substantial data, either online during deployment or offline during training.

In fluid mechanics, adjoint-based MPC has been applied to Navier--Stokes wind-farm models for power maximization through induction and yaw control~\cite{vali2017adjoint, van2022adjoint}, though these works derive problem-specific adjoint equations and do not consider unknown physical parameters. More recently, Alhashim et al.~\cite{alhashim2025} used automatic differentiation through a full Navier--Stokes solver in JAX to optimize swimming gaits and chaotic mixing protocols, while Sunol et al.~\cite{sunol2025learning} implemented differentiable simulation to non-Newtonian fluids by learning constitutive laws from partial flow measurements. The present work builds on this computational infrastructure to couple receding-horizon control with online parameter estimation.

\section{Differentiable Physics for Control and Estimation}
\label{sec:diffphys_control_est}

This section formalizes the control and online estimation framework described in the introduction. We present the formulation in discrete time, with the understanding that in our fluid dynamics example, the ``state'' represents a high-dimensional discretized field (e.g., velocity fields).

\subsection{Discrete-time dynamics with parameters}
Let $x_t \in \mathbb{R}^{n_x}$ denote the state at time index $t \in \{0, 1, 2, \dots\}$, $u_t \in \mathbb{R}^{n_u}$ the control input, and $y_t \in \mathbb{R}^{n_y}$ the measurement. We distinguish
\begin{itemize}
    \item \emph{known} model parameters $\phi \in \mathbb{R}^{n_\phi}$ (e.g., geometry, boundary conditions, known material properties), and
    \item \emph{unknown} physical parameters $\theta \in \mathbb{R}^{n_\theta}$ to be estimated online (e.g., viscosity, constitutive-law parameters, drag coefficients).
\end{itemize}
We assume discrete-time dynamics and measurements:
\begin{align}
x_{t+1} &= f\!\left(x_t, u_t; \phi, \theta\right), 
\label{eq:dyn_general}\\
y_t &= h\!\left(x_t, u_t; \phi, \theta\right).
\label{eq:meas_general}
\end{align}

In many PDE-based models of the system above, the state $x_t$ may be a discretized field over a spatial domain and $y_t$ consists of sparse samples (or a low-dimensional functional) of that field. A common structure is
$
y_t = \mathcal{P}\!\left(x_t\right) ,
%\label{eq:projection_measurement}
$
where $\mathcal{P}$ is a projection/sampling operator (e.g., point sensors, sub-sampled grid, line integrals, camera-based observation model, etc.). We assume the differentiable simulator captures the structural form of $f$ and $h$, while prediction errors arise from mismatch in the unknown parameters $\theta$.

\subsection{Constraints}
We consider standard MPC constraints, including box constraints on inputs and their rates:
\begin{equation}
 u_{\min} \le u_t \le u_{\max},  \quad
-\Delta u_{\max} \le u_t - u_{t-1} \le \Delta u_{\max},
\label{eq:input_constraints}
\end{equation}
and generic inequality constraints
\begin{equation}
g\!\left(x_t, u_t; \phi, \theta\right) \le 0,
\label{eq:general_constraints}
\end{equation}
which may encode actuator limits, safety envelopes, or application-specific constraints (e.g., stress limits, maximum shear rate, bounds on a scalar field, etc.). In partially observed settings, such constraints are typically enforced on predicted/estimated states or outputs.

\subsection{Finite-horizon objective}
Let $\hat x_{t|t}$ denote the state estimate at time $t$ given information up to time $t$, and let $\hat\theta_t$ denote the parameter estimate available at time $t$. For horizon $N$, we define the predicted rollout when $\hat x_{t|t}$ is given
\begin{align}
\hat x_{t+k+1|t} &= f\!\left(\hat x_{t+k|t}, u_{t+k|t}; \phi, \hat\theta_t\right),
\quad k=0,\dots,N-1,
\label{eq:pred_rollout}
\end{align}
with predicted outputs $\hat y_{t+k|t} = h(\hat x_{t+k|t}, u_{t+k|t}; \phi, \hat\theta_t)$.
A general MPC objective can be written as
\begin{multline}
J_t(\mathbf{u}_t;\hat x_{t|t},\hat\theta_t)
=
\sum_{k=0}^{N-1} \ell\!\left(\hat y_{t+k|t}, u_{t+k|t}; r_{t+k}\right)
+ \\ \ell_f\!\left(\hat y_{t+N|t}; r_{t+N}\right),
\label{eq:mpc_cost_general}
\end{multline}
where $\mathbf{u}_t := \{u_{t|t},\dots,u_{t+N-1|t}\}$ and $r_{t+k}$ denotes a reference/command (or desired field statistic).

\subsection{Simulation-based MPC via differentiable physics}
At each time $t$, the MPC problem is
\begin{subequations} \label{eq:mpc_problem}
\begin{align}
\min_{\mathbf{u}_t} \quad &
J_t(\mathbf{u}_t;\hat x_{t|t},\hat\theta_t)
\label{eq:mpc_cost}\\
\text{s.t.}\quad &
\hat x_{t+k+1|t} = f\!\left(\hat x_{t+k|t}, u_{t+k|t}; \phi, \hat\theta_t\right), \label{eq:mpc_dync}\\
& \hat y_{t+k|t} = h\!\left(\hat x_{t+k|t}, u_{t+k|t}; \phi, \hat\theta_t\right),
\label{eq:mpc_output_model}\\
&
u_{\min} \le u_{t+k|t} \le u_{\max}, \label{eq:mpc_input_bounds}\\
-&\Delta u_{\max} \le u_{t+k|t} - u_{t+k-1|t} \le \Delta u_{\max},
\label{eq:mpc_delta_input_bounds}\\
&
g\!\left(\hat x_{t+k|t}, u_{t+k|t}; \phi, \hat\theta_t\right) \le 0, \label{eq:mpc_generic_constraint}
\end{align}
\end{subequations}
where $ k=0,\dots,N-1$. Here $u_{t-1}$ denotes the control input most recently applied to the plant. For $k=0$, the rate constraint is understood with $u_{t-1|t}:=u_{t-1}$. Let $\mathbf u_t^\star := \{u_{t|t}^\star,\dots,u_{t+N-1|t}^\star\}$ denote the optimizer of \eqref{eq:mpc_problem}; the control applied to the plant is the first element, $u_t = u_{t|t}^\star$.

\begin{figure*}[t]
\centering
\begin{tikzpicture}[
    >=Latex,
    node distance=1.4cm and 1.5cm,
    estblock/.style={
        draw, rounded corners, thick, align=center,
        text width=5.16cm, minimum height=2.35cm, fill=white,
        inner xsep=8pt, inner ysep=4.5pt
    },
    mpcblock/.style={
        draw, rounded corners, thick, align=center,
        text width=4.35cm, minimum height=1.9cm, fill=white,
        inner sep=4.5pt
    },
    plantblock/.style={
        draw, rounded corners, thick, align=center,
        text width=3.55cm, minimum height=1.7cm, fill=white,
        inner sep=4.5pt
    },
    line/.style={-Latex, thick},
    every node/.style={font=\small}
]

% Blocks
\node[estblock] (mhe) {%
\textbf{State estimation}~\eqref{eq:mhe_problem}\\[0.8mm]
$\displaystyle
\begin{aligned}
\min_{\mathbf{z}_{t-L:t}} &
\sum_{k=t-L}^{t}
\|y_k-h(z_k,u_k;\phi,\hat\theta_t)\|_{R}^{2}\\
&\qquad+ \|z_{t-L}-\bar x_{t-L}\|_{P}^{2}
\end{aligned}
$\\[0.8mm]
subject to \eqref{eq:mhe_constraints}
};

\node[mpcblock, right=1.7cm of mhe] (mpc) {%
\textbf{MPC}~\eqref{eq:mpc_problem}\\[0.8mm]
$\displaystyle
\mathbf{u}_{t}^{\star}
=
\arg\min_{\mathbf{u}_t}
J_t(\mathbf{u}_t;\hat x_{t|t},\hat\theta_t)
$\\[0.8mm]
subject to \eqref{eq:mpc_dync}--\eqref{eq:mpc_generic_constraint}
};

\node[plantblock, right=1.6cm of mpc] (plant) {%
\textbf{Physical plant}\\[0.8mm]
$\displaystyle x_{k+1}=f(x_k,u_k;\phi,\theta)$\\
$\displaystyle y_k=h(x_k,u_k;\phi,\theta)$
};

\node[estblock, below=0.5cm of mpc] (param) {%
\textbf{Parameter estimation}~\eqref{eq:param_est}\\[0.8mm]
$\displaystyle
\begin{aligned}
\min_{\theta\in\Theta} &
\sum_{k=t-W}^{t}
\|y_k - \hat y_{k|t}\|_{R_\theta}^{2} + \lambda_\theta\|\theta-\bar\theta\|^{2}
\end{aligned}
$\\[0.8mm]
subject to \eqref{eq:param_est_constraint}-\eqref{eq:param_est_output}
};

% Main horizontal lines
\draw[line] (mhe.east) -- node[above] {$\hat x_{t|t}$} (mpc.west);
\draw[line] (mpc.east) -- node[above] {$u_t = u_{t|t}^{\star}$} (plant.west);

% Plant to parameter estimation
\draw[line] (plant.south) |- node[pos=0.25,right] {$y_t$} (param.east);

% State estimation to parameter estimation
\draw[line]
    ($(mhe.south)+(0.15,0)$)
    |- node[pos=0.24,left] {$\hat x_{t-L:t}$}
    ($(param.west)+(0,-0.32)$);

% Parameter estimation to state estimation
\coordinate (thetastart) at ($(param.west)+(0,0.32)$);
\coordinate (thetabranch) at ($(thetastart)+(-0.5,0)$);

\draw[line]
    (thetastart) -- (thetabranch)
    -| node[pos=0.24,above] {$\hat\theta_t$}
    ($(mhe.south)+(0.95,0)$);

% Branch of parameter estimate to MPC
\draw[line]
    (thetabranch)
    |- node[pos=0.32,right] {}
    ($(mpc.west)+(0,-0.55)$);

\end{tikzpicture}
\caption{Control and estimation with a shared differentiable model, with rollout constraints enforced in each of the three optimizations.}
\label{fig:joint_control_estimation_diagram}
\end{figure*}
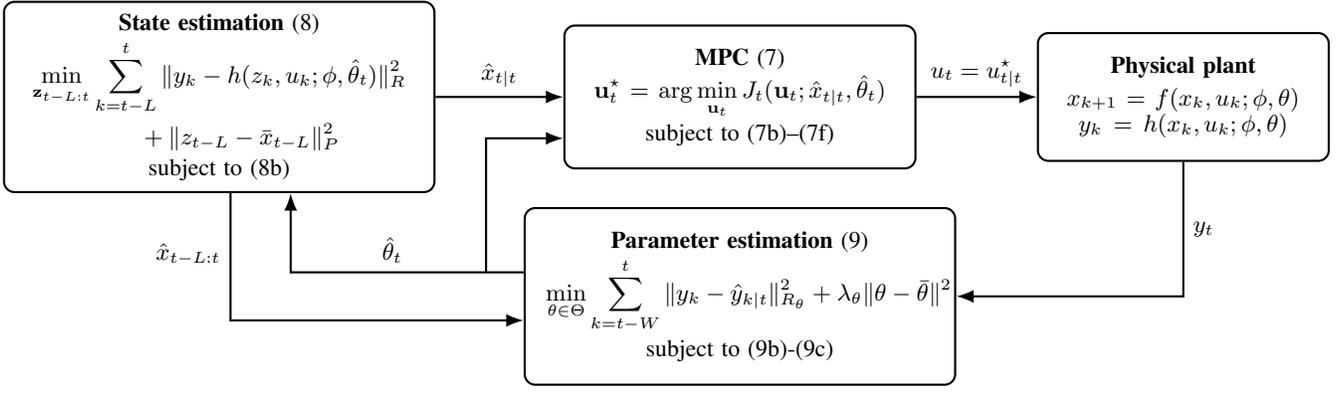

The key distinction from classical MPC is that $f$ and $h$ may not be available in closed algebraic form; instead, they are implemented by a differentiable simulator. We therefore treat the dynamics in \eqref{eq:mpc_dync} as an implicit computational constraint enforced through a rollout \eqref{eq:pred_rollout}, and we compute gradients $\nabla_{\mathbf{u}_t} J_t$ and (optionally) $\nabla_{\hat x_{t|t}} J_t,\ \nabla_{\hat\theta_t} J_t$ by differentiating through the simulation via automatic differentiation.
This enables standard gradient-based solvers without requiring explicit analytic derivatives of the underlying physics. In PDE-based systems, even a modest prediction horizon may correspond to thousands of simulator timesteps, making gradient computation through the full rollout the central computational challenge.

% \textbf{Hyperparameters:} The key hyperparameters in designing \eqref{eq:mpc_problem} are: the time horizon $N$, terms in the cost function (which may include regularization for smoothing of the control), weights in the cost function, and the key constraints. There are hyperparameters in the solver which also must be chosen such as number of epochs, learning rate, schedule, stopping criteria, etc.

\subsection{State estimation under sparse measurements}
When the state is not directly observed, MPC requires a state estimate $\hat x_{t|t}$. For high-dimensional fields and sparse measurements, an optimization-based observer is natural. Using past inputs $\mathbf{u}_{t-L:t-1}:=\{u_{t-L},\dots,u_{t-1}\}$ and measurements $\mathbf{y}_{t-L:t}:=\{y_{t-L},\dots,y_t\}$ over a window of length $L$, and holding $\hat\theta_t$ fixed, we consider a moving-horizon observer that optimizes over a candidate trajectory $\mathbf{z}_{t-L:t}:=\{z_{t-L},\dots,z_t\}$:
\begin{subequations} \label{eq:mhe_problem}
\begin{align}
\min_{\mathbf{z}_{t-L:t}} &
\sum_{k=t-L}^{t} \|y_k - h(z_k, u_k; \phi, \hat\theta_t)\|_{R}^2
+ \|z_{t-L} - \bar x_{t-L}\|_{P}^2
\label{eq:mhe_cost}\\
\text{s.t.}\quad &
z_{k+1} = f(z_k, u_k; \phi, \hat\theta_t),
\quad k=t-L,\dots,t-1.
\label{eq:mhe_constraints}
\end{align}
\end{subequations}
Let $\mathbf{z}_{t-L:t}^\star:=\{z_{t-L}^\star,\dots,z_t^\star\}$ denote the optimizer of \eqref{eq:mhe_problem}. The observer output supplied to MPC is then the terminal estimate $\hat x_{t|t}:=z_t^\star$. Here $\bar x_{t-L}$ is an arrival prior at the beginning of the estimation window, and $\|z\|_M^2 := z^\top M z$.

\subsection{Online parameter estimation with differentiable physics}
To identify unknown physical parameters $\theta$, we use past input/output data over a window $W$, with $\mathbf{u}_{t-W:t-1}:=\{u_{t-W},\dots,u_{t-1}\}$ and $\mathbf{y}_{t-W:t}:=\{y_{t-W},\dots,y_t\}$. Starting from an estimated state at the beginning of the window, denoted $\hat x_{t-W|t}$, we simulate forward and fit $\theta$ by minimizing the mismatch between measured outputs $y_k$ and predicted outputs $\hat y_{k|t}$:
\begin{subequations} \label{eq:param_est}
\begin{align}
\min_{\theta \in \Theta} \quad &
\sum_{k=t-W}^{t} \|y_k - \hat y_{k|t}\|_{R_\theta}^2
+ \lambda_\theta \|\theta - \bar\theta\|^2
\label{eq:param_est_cost}\\
\text{s.t.}\quad &
\hat x_{k+1|t} = f(\hat x_{k|t}, u_k; \phi, \theta),
\quad k=t-W,\dots,t-1,
\label{eq:param_est_constraint}\\
&
\hat y_{k|t} = h(\hat x_{k|t}, u_k; \phi, \theta),
\quad k=t-W,\dots,t.
\label{eq:param_est_output}
\end{align}
\end{subequations}
Here $\Theta$ encodes known bounds or priors, $\bar\theta$ is a prior, and $\lambda_\theta \ge 0$ regularizes drift. Let $\theta^\star$ denote the optimizer of \eqref{eq:param_est}; the parameter estimate supplied to MPC is then $\hat\theta_t := \theta^\star$.

\subsection{Coupling MPC, state estimates, and parameter estimates}
The MPC problem \eqref{eq:mpc_problem} is solved at each time $t$ using the current state and parameter estimates. In partially observed systems, this requires a choice of estimation architecture: a \emph{fully coupled} scheme updates both $\hat x_{t|t}$ and $\hat\theta_t$ at every step, while a \emph{decoupled} scheme computes $\hat x_{t|t}$ synchronously (e.g., via \eqref{eq:mhe_problem}) and updates $\hat\theta_t$ asynchronously every $M$ steps by solving \eqref{eq:param_est} over a longer window. Under the decoupled architecture, the loop is: estimate $\hat x_{t|t}$, solve \eqref{eq:mpc_problem} using $(\hat x_{t|t},\hat\theta_t)$, apply $u_t=u_{t|t}^\star$, then update $\hat\theta_t$ periodically.

In the accelerator example, the behavior is modeled as a static map, so state estimation is not required and $\hat\theta_t$ is updated synchronously at each control step. In the fluid mechanics example, we likewise assume the state is known. A key challenge is that closed-loop control may suppress informative transients, causing the parameter estimator in \eqref{eq:param_est} to stagnate under sparse measurements. We do not study observability here; in practice, window design, regularization, and excitation terms in \eqref{eq:mpc_cost_general} can help mitigate this issue.

\section{Fluid Mechanics Example} \label{sec:fluid_example}
\subsection{Problem setup}

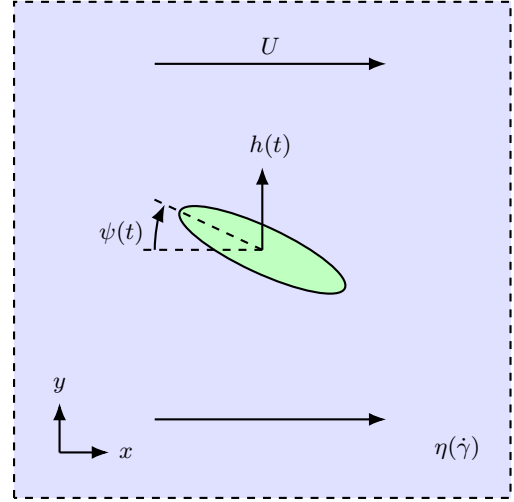
\begin{figure}[t]
\centering
\begin{tikzpicture}[
    >=Latex,
    scale=1.1,
    every node/.style={font=\small}
]

% Domain
\fill[blue!12] (0,0) rectangle (6,6);
\draw[thick,dashed] (0,0) rectangle (6,6);

% Flow arrows
\draw[->, thick] (1.7,5.25) -- (4.5,5.25) node[midway, above] {$U$};
\draw[->, thick] (1.7,0.95) -- (4.5,0.95);

% Viscosity label
\node[anchor=south east] at (5.75,0.35) {$\eta(\dot{\gamma})$};

% Coordinate axes
\draw[->, thick] (0.55,0.55) -- (1.15,0.55) node[right] {$x$};
\draw[->, thick] (0.55,0.55) -- (0.55,1.15) node[above] {$y$};

% Swimmer center
\coordinate (C) at (3,3);

% Swimmer body
\begin{scope}[rotate around={-25:(C)}]
    \fill[green!25] (C) ellipse [x radius=1.1, y radius=0.28];
    \draw[thick] (C) ellipse [x radius=1.1, y radius=0.28];
\end{scope}

% Heave arrow
\draw[->, thick] (C) -- ++(0,1.0);
\node at ($(C)+(0.10,1.25)$) {$h(t)$};

% Dashed angle reference lines
\draw[dashed, thick] (C) -- ++(-1.45,0);
\draw[dashed, thick] (C) -- ++({1.45*cos(155)},{1.45*sin(155)});

% Angle arc and label
\draw[->, thick] ($(C)+(-1.3,0)$) arc[start angle=180,end angle=155,radius=1.3];
\node at ($(C)+(-1.7,0.2)$) {$\psi(t)$};

\end{tikzpicture}
\caption{The elliptical foil in a 
background flow $U$, with heave $h(t)$ and pitch $\psi(t)$ parameterized by the Fourier series in \eqref{eq:heave} and \eqref{eq:pitch}. The solver implements the generalized Newtonian constitutive model \eqref{eq:carreau_yasuda}, though results in this work use the Newtonian limit ($\eta_0 = \eta_\infty = \eta$). The domain has periodic boundary conditions.}
\label{fig:schematic}
\end{figure}

We consider a two-dimensional elliptical foil (semi-axes $a=0.5$, 
$b=0.06$, chord $c=2a=1.0$) holding station in a
background flow, as illustrated in Fig.~\ref{fig:schematic}.
A mean background flow of velocity $U$ in the $+x$ direction is maintained
by a pressure correction applied at each timestep, and the swimmer
controls its vertical displacement (heave) and angular orientation (pitch)
to hold station efficiently. The oscillatory body velocity 
$\mathbf{v}_{\text{body}}(t)$ arising from the heave and pitch motions 
is prescribed by the gait parameterization. The simulation is warmed up for 8 flapping periods before 
the MPC loop begins, ensuring a developed wake at the start of replanning.

Following~\cite{alhashim2025}, the heaving and pitching motions are
parameterized by truncated Fourier series:
\begin{align}
h(t) &= \sum_{m=1}^{M_f} \alpha_m \cos(2\pi m f t) + \beta_m \sin(2\pi m f t), \label{eq:heave}\\
\psi(t) &= \psi_0 + \sum_{m=1}^{M_f} a_m \cos(2\pi m f t + \zeta) \nonumber \\
&\qquad\qquad + b_m \sin(2\pi m f t + \zeta), \label{eq:pitch}
\end{align}
where $f$ is the base frequency, $\zeta$ is a heave--pitch phase lag,
and $\psi_0$ is a constant pitch offset. The $m=0$ (DC) modes are
excluded from optimization: a constant heave offset is redundant with the
particle position, and a constant pitch offset can lead to degenerate
static-airfoil solutions that reduce drag without producing thrust. The variable $t$ in this section is the continuous time variable, not to be confused with the discrete index in the previous section.

The fluid is governed by the 2D incompressible momentum and continuity equations,
\begin{align}
\rho\left(\frac{\partial \mathbf{v}}{\partial t} + \mathbf{v}\cdot\nabla\mathbf{v}\right)
&= -\nabla p + \nabla\cdot\boldsymbol{\tau} + \mathbf{f}_{\text{IBM}}, \label{eq:momentum}
\\ 
\nabla\cdot\mathbf{v}&= 0,
\label{eq:continuity}
\end{align}
where $\mathbf{v} = [v_x(x,y), v_y(x,y)]^\top$ is the velocity field, $p(x,y)$ is the pressure, and
$\mathbf{f}_{\text{IBM}}(x,y)$ is the immersed boundary forcing that
enforces the no-slip condition on the swimmer surface. The fluid density is $\rho=1$ and the gradient operator is given by $\nabla \coloneqq [\partial_x, \partial_y]^\top$ --- throughout this section, $x,y$ denotes the independent spatial variables of the PDE, and should be distinguished from the state and measurement vectors $x_k, y_k$ in \eqref{eq:dyn_general}-\eqref{eq:meas_general} which appear with subscripts. The deviatoric
stress tensor is $\boldsymbol{\tau} = \eta(\dot\gamma)\left(\nabla\mathbf{v} + \nabla\mathbf{v}^\top\right)$,
where the generalized viscosity $\eta$ depends on the local shear
rate $\dot\gamma = \sqrt{2\, \text{tr}(\mathbf{D} \mathbf{D}^\top)}$, and where $\mathbf{D} = \tfrac{1}{2}(\nabla\mathbf{v}+\nabla\mathbf{v}^\top)$
is the rate-of-strain tensor. For a Newtonian fluid, $\eta$ is a
constant; for a generalized Newtonian fluid, it is given by a
constitutive model such as the Carreau--Yasuda relation:
\begin{equation}
\eta(\dot\gamma) = \eta_\infty + (\eta_0 - \eta_\infty)\left[1 + (\lambda\dot\gamma)^a\right]^{(n-1)/a},
\label{eq:carreau_yasuda}
\end{equation}
where $\eta_0$ and $\eta_\infty$ are the zero- and
infinite-shear-rate viscosities, $\lambda$ is a relaxation time,
$n$ is the power-law index, and $a$ controls the transition
sharpness. The Newtonian limit is recovered when $\eta_0 = \eta_\infty = \eta$,
in which case $\boldsymbol{\tau}$ reduces to the Newtonian stress tensor
\eqref{eq:continuity} reduce to the
Navier--Stokes equations.

Mapping to the general framework of Section~\ref{sec:diffphys_control_est}:
the state vector $x_k$ in \eqref{eq:dyn_general} is a collection of the discretized velocity field $\mathbf{v}$ on an $256 \times 256$ grid of points on spatial variables $(x,y)$ at time $k$. The control inputs $u_k$ are the Fourier gait coefficients
$\{a_m, b_m, \alpha_m, \beta_m\}_{m=1}^{M_f}$. These coefficients are held fixed over the horizon, so that the entire input trajectory is generated by a single periodic gait. From the perspective of the general MPC formulation \eqref{eq:mpc_problem}, this restriction to Fourier-parameterized, horizon-constant gaits can be interpreted as an application-specific constraint on the admissible inputs, consistent with the generic constraint class in \eqref{eq:mpc_generic_constraint}. The known parameters $\phi$ include the domain geometry, grid resolution,
boundary conditions, and base frequency; and the unknown parameter
$\theta$ is the rheological quantity estimated online---in the
Newtonian case, the scalar viscosity $\eta$.

The measurement available to the swimmer consists of the hydrodynamic
forces on its body, computed via the immersed boundary method (IBM).
At each timestep, the IBM computes the force density $\mathbf{f}_{\text{IBM}}$
in~\eqref{eq:momentum} required to enforce the no-slip condition on the
swimmer surface. The net hydrodynamic force on the swimmer is
$\mathbf{F} = [F_x, F_y]^\top = \oint \mathbf{f}_{\text{IBM}} \, dS$.
From this, we define two scalar measurements per flapping period $T=1/f$: the mean thrust power $\bar P_{\text{thrust}} = U \int_0^T F_x \, dt / T$
and the mean oscillatory input power
$\bar P_{\text{osc}} = \frac{1}{T}\int_0^T \!\oint \mathbf{f}_{\text{IBM}} \cdot \mathbf{v}_{\text{body}} \, dS \, dt$,
where $\mathbf{v}_{\text{body}}$ is the prescribed body velocity. These are the quantities a physical
swimmer could sense via embedded strain or force sensors.

The gait optimization objective is to maximize the Froude-type swimming
efficiency~\cite{alhashim2025}:
\begin{equation}
\mathcal{E} = \bar P_{\text{thrust}}/\bar P_{\text{osc}},
\label{eq:froude_efficiency}
\end{equation}
by adjusting the Fourier coefficients, while simultaneously estimating
$\theta$ from the mismatch between predicted and observed power
measurements. In our implementation, the MPC problem \eqref{eq:mpc_problem}
is solved over a horizon of $4{,}000$ simulation time steps, corresponding to two flapping periods. The objective is evaluated by
computing the average efficiency \eqref{eq:froude_efficiency} over only
the final $2{,}000$ steps of the rollout, so that the optimization
targets a steady-state gait after the transient dynamics in the first
period have decayed.

Since the full simulator state $\hat x_{t|t} = x_{t}$ is available at each replanning step, explicit state estimation (e.g., as in \eqref{eq:mhe_problem}) is therefore not required, and in future work we will consider the harder problem where the full state must be estimated from sparse field measurements.

\subsection{Implementation}

The plant model is a JAX-based incompressible flow solver
with immersed boundary forcing and a generalized Newtonian viscosity model \cite{sunol2025learning}, although only the Newtonian case is used in the present demonstration. The entire forward model---including
the IBM force computation, pressure projection, and constitutive
evaluation---is end-to-end differentiable, providing exact control
gradients $\nabla_u \mathcal{E}$ for gait optimization and parameter
sensitivities $\nabla_\theta$ of the estimation 
objective~\eqref{eq:param_est} via automatic differentiation.

The receding-horizon loop follows the asynchronous architecture of
Section~\ref{sec:diffphys_control_est}. At each replanning step, a
candidate policy (Fourier gait) is evaluated over a prediction horizon
of two flapping periods ($4{,}000$ timesteps). The cost~\eqref{eq:froude_efficiency} is computed only over the second period, after the flow has reached
steady state under the candidate policy; gradients are propagated through the full two-period horizon, including 
the transient first period. The policy is then updated via the
Adam optimizer~\cite{kingma2014adam}, and both periods are committed to 
the plant: the first serves as a wake-development warmup under the new 
gait, and the second is the measurement period from which power data is 
recorded for parameter estimation. The controller then replans from the 
resulting state, so each MPC cycle advances the system by two flapping 
periods.  At each gait switch, the incoming Fourier coefficients are 
phase-matched to the current swimmer state so that $h(t)$ and 
$\psi(t)$ remain continuous across the transition. The restriction 
of the loss to the second period is deliberate: optimizing over the 
first (transient) period would reward gaits that exploit wake startup 
artifacts rather than steady-state efficiency.

Parameter estimation is performed asynchronously every $M=5$ control
updates. A sliding buffer of size $W=5$ stores the plant state, applied 
policy, and measured power output from each committed two-period cycle. 
The estimator minimizes the power-prediction 
mismatch~\eqref{eq:param_est} over this buffer: each buffered cycle is 
re-simulated under the candidate $\hat\theta$, and the predicted power 
is compared to the recorded measurement. With $5$ buffered cycles of 
$4{,}000$ timesteps each, the estimator differentiates through 
$20{,}000$ total timesteps per update.
For the Newtonian case ($\theta = \eta$), we
reparameterize as $\log\eta$ to maintain positivity without box
constraints. All simulations use a timestep
$\Delta t = 10^{-3}$.

The gait optimizer runs 5 Adam steps per MPC cycle with learning 
rates of $1.5$--$2 \times 10^{-2}$ (rotation) and 
 $3$--$4 \times 10^{-3}$ (displacement); the viscosity estimator 
 runs 1 gradient step per update with learning rates ranging from 
 10 to 40 across configurations. The three configurations in 
Fig.~\ref{fig:results} differ in gait optimizer variant (Adam 
 with and without backtracking line search), viscosity learning 
 rate, and whether a late-stage learning rate decay (factor 0.25 
 after 20 MPC cycles) is applied. These choices were not 
 extensively tuned.
 
\subsection{Results}

\begin{figure}[t]
\centering
\includegraphics[width=0.9\columnwidth]{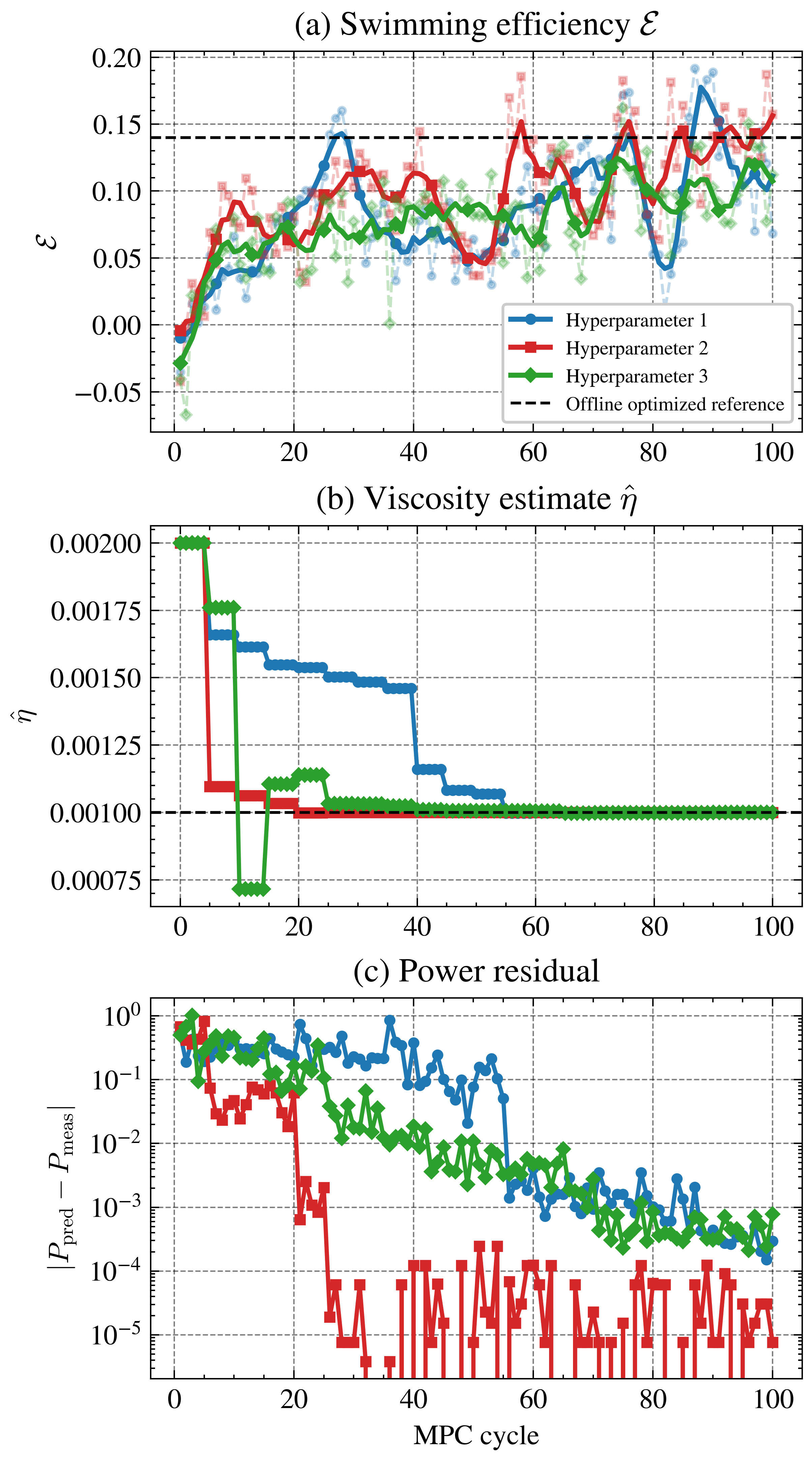}
\caption{Joint MPC gait optimization and viscosity estimation at $Re = 1{,}000$ over 100 replanning cycles,
shown for three hyperparameter configurations that differ in learning rates and late-stage decay schedule.
(a)~Swimming efficiency $\mathcal{E}$; solid lines are
5-point centered moving averages and translucent markers show
the raw per-cycle values; the dashed line is an offline 
gait optimization with known $\eta$ over many gradient steps.
(b)~Viscosity estimate $\hat\eta$; the dashed line is
$\eta_{\text{true}} = 10^{-3}$.
(c)~Absolute power-prediction residual
$|P_{\text{pred}} - P_{\text{meas}}|$.}
\label{fig:results}
\end{figure}

We demonstrate the joint control and estimation framework in the
Newtonian limit, estimating the scalar viscosity $\eta$ while
simultaneously optimizing the swimming gait. The true viscosity is
set to $\eta_{\text{true}} = 10^{-3}$, corresponding to
$Re = Uc/\eta = $ $1{,}000$ and the estimator is initialized at
$\hat\eta_0 = 2\eta_{\text{true}}$. The swimmer begins with a
single-frequency heave--pitch gait and $M_f = 6$ Fourier modes
available for optimization.

Figure~\ref{fig:results} shows the evolution of the system over
$100$ MPC cycles. The viscosity estimate $\hat\eta$
converges to within $1\%$ of the true value within
$4$--$10$ control updates depending on hyperparameter configuration (Fig.~\ref{fig:results}b), at which
point the gait optimizer has access to an accurate plant model. The power-prediction residual $|\bar P_{\text{pred}}(\hat\eta) - \bar P_{\text{meas}}|$ falls below $10^{-3}$ for all three configurations (Fig.~\ref{fig:results}c).

Two features of the coupled dynamics are worth noting. First, the
early MPC cycles produce policies that are suboptimal not because
the optimizer fails, but because it solves the \emph{wrong} problem:
with $\hat\eta > \eta_{\text{true}}$, the predicted flow is more
viscous than reality, and the optimizer designs gaits suited to a
lower-$Re$. Despite this model mismatch, the optimizer still improves
efficiency during the early cycles because the qualitative physics
of the correct Reynolds-number regime are preserved. However,
efficiency oscillations persist well after $\hat\eta$ has converged
(Fig.~\ref{fig:results}a), indicating that the residual noise is
driven by the gait optimization itself---nonconvexity of the
landscape and aggressive per-step updates---rather than by
parameter estimation error.

Second, the gait changes driven by the MPC loop naturally provide
persistence of excitation for the parameter estimator. Each
replanning step alters the velocity and shear-rate distribution
around the swimmer, probing different regions of the constitutive
response. This is in contrast to a fixed-gait scenario, where the
estimator would observe the same flow conditions repeatedly and
convergence of $\hat\eta$ would depend entirely on the information
content of a single operating point. The receding-horizon structure
thus provides a self-exciting property: control updates that improve
$\mathcal{E}$ simultaneously enrich the data available for parameter
estimation.

The MPC-optimized gaits achieve swimming efficiencies that
approach the offline-optimized reference of
$\mathcal{E} \approx 0.14$, up from near-zero or negative
initial values (negative $\mathcal{E}$ indicates the body produces 
net drag rather than thrust) under the starting gait and incorrect viscosity
estimate. The persistent efficiency oscillations visible in
Fig.~\ref{fig:results}(a) reflect several additional sources of
approximation beyond gait nonconvexity: the short single-period
replanning horizon, continuity corrections
applied to the Fourier coefficients at each gait switch to
prevent discontinuities in the swimmer trajectory, and the
deliberate choice not to extensively tune the optimizer
hyperparameters for this demonstration. The framework extends
directly to Carreau--Yasuda parameters; extending it to partial
state observations and simultaneous estimation of
$(\eta_0, \lambda, n)$ is deferred to future work.

\section{Accelerator System: A Static Example} \label{sec:acc_example}
\subsection{Static Formulation}
Particle accelerator tuning is a useful case for our framework, since high-fidelity models are generally computational and are rarely used directly in control. We reformulate the setup of Section~\ref{sec:diffphys_control_est} for this setting because, at the controller update rate, the plant appears \emph{static}. In such tuning problems, a particle bunch traverses the beamline essentially instantaneously relative to actuator updates: the controller sets magnet setpoints, the bunch propagates through the beamline, and diagnostics return measurements. Only after the bunch has propagated and the measurements have been made can the next control action be applied. The accelerator is therefore naturally modeled as a static input--output map.

Accordingly, from \eqref{eq:dyn_general}--\eqref{eq:meas_general} we collapse the fast internal evolution into the measurement model and write
\begin{equation}
y_t = h\!\left(u_t; \phi, \theta\right),
\label{eq:acc_static_map}
\end{equation}
where $u_t\in\mathbb{R}^{n_u}$ are the setpoints held during bunch $t$ and $y_t\in\mathbb{R}^{n_y}$ are the resulting diagnostic readbacks. Unmeasured initial beam conditions are absorbed into the parameters: $\phi$ denotes known configuration quantities, while $\theta$ captures unknown or drifting effects to be estimated online. The differentiable simulator then serves as a forward map from $u_t$ to predicted measurements $\hat y_t$.

\subsection{MPC in the static setting}
Although \eqref{eq:acc_static_map} is static in time $t$, an MPC formulation remains valuable. Many accelerator objectives are sequential, for example enforcing smooth setpoint moves and rate limits over multiple steps. Likewise, safety constraints are often expressed in diagnostic space and can be imposed directly on predicted measurements. Planning over a horizon therefore provides a look-ahead mechanism for navigating a generally nonconvex tuning problem that depends on the current parameter estimate $\hat\theta_t$.

Given a horizon $N$ and a current parameter estimate $\hat\theta_t$, we define predicted diagnostic outputs
\begin{equation}
\hat y_{t+k|t} = h\!\left(u_{t+k|t}; \phi, \hat\theta_t\right),
\qquad k=0,\dots,N.
\label{eq:acc_pred_outputs}
\end{equation}
A generic objective consistent with \eqref{eq:mpc_cost_general} becomes
\begin{equation}
J_t(\mathbf{u}_t;\hat\theta_t)
=
\sum_{k=0}^{N-1}
\ell\!\left(\hat y_{t+k|t},u_{t+k|t}\right)
+ \ell_f\!\left(\hat y_{t+N|t}\right),
\label{eq:acc_cost}
\end{equation}
where $\mathbf{u}_t := \{u_{t|t},\dots,u_{t+N-1|t}\}$. Application-specific constraints are imposed directly on the predicted diagnostics and inputs through a differentiable inequality map
\begin{equation}
g\!\left(\hat y_{t+k|t},u_{t+k|t};\phi,\hat\theta_t\right) \le 0,
\qquad k=0,\dots,N-1.
\label{eq:acc_output_constraints}
\end{equation}

Putting these together, the static-horizon MPC problem at time $t$ is
\begin{subequations}\label{eq:acc_static_mpc}
\begin{align}
\min_{\mathbf{u}_t}\quad & J_t(\mathbf{u}_t;\hat\theta_t) \label{eq:acc_static_mpc_cost}\\
\text{s.t.}\quad &
\hat y_{t+k|t} = h\!\left(u_{t+k|t};\phi,\hat\theta_t\right), \label{eq:acc_static_mpc_model}\\
& u_{\min} \le u_{t+k|t} \le u_{\max}, \label{eq:acc_static_mpc_bounds}\\
& -\Delta u_{\max} \le u_{t+k|t} - u_{t+k-1|t} \le \Delta u_{\max}, \label{eq:acc_static_mpc_rate}\\
& g\!\left(\hat y_{t+k|t},u_{t+k|t};\phi,\hat\theta_t\right) \le 0, \label{eq:acc_static_mpc_g}
\end{align}
\end{subequations}
for $k=0,\dots,N-1$. Here $u_{t-1}$ denotes the most recently applied setpoint, and for $k=0$ we define $u_{t-1|t}:=u_{t-1}$. With slight abuse of notation, $\mathbf{u}_t$ denotes the candidate horizon input sequence; the optimized sequence is denoted by $\mathbf{u}_t^\star$. Because $h$ is provided by a differentiable beamline simulator, gradients of the objective and constraints with respect to the horizon inputs are obtained by differentiating through the simulation.

\subsection{Online parameter estimation in the static setting}
The same static viewpoint naturally yields a simplified parameter identification problem. Using a window of past applied setpoints $\{u_k\}$ with measured diagnostics $\{y_k\}$, we update $\theta$ by minimizing a simulation-to-measurement mismatch:
\begin{equation}
\min_{\theta\in\Theta}\quad
\sum_{k=t-W}^{t}
\left\|y_k - h\!\left(u_k;\phi,\theta\right)\right\|_2^2
+ \lambda_\theta \|\theta-\bar\theta\|^2,
\label{eq:acc_param_est}
\end{equation}
where $\Theta$ encodes known parameter bounds/priors and $\bar\theta$ is a nominal parameter vector, and the optimizer of \eqref{eq:acc_param_est} is denoted as $\hat \theta_t$. As in Sec.~\ref{sec:diffphys_control_est}, differentiability of $h$ provides $\nabla_\theta$ efficiently, enabling frequent, warm-started updates that track calibration drift and model mismatch during operation.

\subsection{Problem Setup}
Our differentiable model is built using Cheetah \cite{kaiser2024cheetah, stein2022accelerating}, which is a GPU accelerated particle beam modeling code. It is able to model general beamlines and contains 20+ unique beamline elements, including models for challenging space charge effects, with auto-differentiability provided by the PyTorch backend. The beamline is modeled after the Isotope Production Facility (IPF) beamline, which is an approximately 30 meter long beamline containing 9 quadrupole magnets for beam focusing and 4 bending magnets \cite{osti_977750} to direct the beam off of the primary linear accelerator axis. The IPF beamline extracts 100 MeV protons from the primary linear accelerator and directs them to a target, which produces radioisotopes for medical research, diagnosis, and other uses. It is essential that the shape of the particle beam is controlled to irradiate the surface of the target evenly, while also ensuring that the size of the particle bunch does not exceed the dimensions of the beam pipe along the way to the target in order to avoid damage to accelerator components.

The control input is the vector of the first 7 quadrupole strengths
$u_t \in \mathbb{R}^{7}$,
corresponding to the field strength setpoints of each quadrupole
$\text{Q1},\dots,\text{Q7}$.
The measurement vector $y_t \in \mathbb{R}^{8}$ consists of transverse beam sizes at four diagnostic locations (segment exits) along the beamline:
\begin{equation}
y_t =
\big[
\sigma_x^{a}, \sigma_y^{a},
\sigma_x^{b}, \sigma_y^{b},
\sigma_x^{c}, \sigma_y^{c},
\sigma_x^{d}, \sigma_y^{d}
\big]^\top,
\label{eq:ipf_meas_vec}
\end{equation}
where $\sigma_{x/y}^{(\cdot)}$ are the RMS sizes (in the horizontal and vertical directions) at the end of segments $a$--$d$ within the 30 meter length.
The unknown parameter vector $\theta\in\mathbb{R}^4$ parametrizes the incoming
beam distribution and is taken as
$
    \theta = [\sigma_x,\sigma_{p_x},\sigma_y,\sigma_{p_y}]^\top
$
and is generally unknown in practice. The measurements do not exactly correspond to the current physical placement of these sensors at IPF, as this is a preliminary simulation based study.

We use a tracking loss on the diagnostic nearest the target to achieve a circular 10 mm beam:
\begin{multline}   
J_t(\mathbf{u}_t;\hat\theta_t)
=
\frac{1}{N}\sum_{k=0}^{N-1}
\left\|
\hat y_{t+k|t}^{(d)} - y^\star
\right\|_2^2 \\
\;+\; 
\lambda_{\mathrm{smooth}}\,\Omega_{\mathrm{smooth}}
\;+\;
\lambda_{\mathrm{soft}}\,\Omega_{\mathrm{soft}},
\label{eq:ipf_cost_used}
\end{multline}
where $\hat y_{t+k|t}^{(d)} = [\hat\sigma_x^{d},\hat\sigma_y^{d}]^\top$
denotes the predicted beam sizes at the final diagnostic location,
and $y^\star = [\sigma^\star,\sigma^\star]^\top$ with $\sigma^\star = 10~\mathrm{mm}$. We emphasize that the stage term tracks the \emph{terminal diagnostic location} $d$ at every predicted step within the horizon, rather than tracking all intermediate diagnostics. The following loss terms penalize rate changes to the planned control inputs (quadrupole setpoints) and predicted beam size measurements greater than the threshold:
\begin{align}
\Omega_{\mathrm{smooth}}
&=
\left\|u_{t|t}-u_{t-1}\right\|_2^2
+
\sum_{k=1}^{N-1}\left\|u_{t+k|t}-u_{t+k-1|t}\right\|_2^2,
\label{eq:ipf_smooth} \\
\Omega_{\mathrm{soft}}
&=
\sum_{k=0}^{N-1}
\left\|
\max\!\left\{0,\;\hat y_{t+k|t}-y_{\max}\right\}
\right\|_2^2,
\label{eq:ipf_soft_constraint}
\end{align}
where $N = 5$, $\lambda_{\mathrm{smooth}} = 10^{-4}$, $\lambda_{\mathrm{soft}}=1.0$, and threshold $y_{\max}=14$~mm. Here $u_{t-1}$ denotes the setpoint applied on the most recent step, i.e., the current actuator state at time $t$.

\begin{figure}[t]
\centering
\includegraphics[width=0.9\columnwidth]{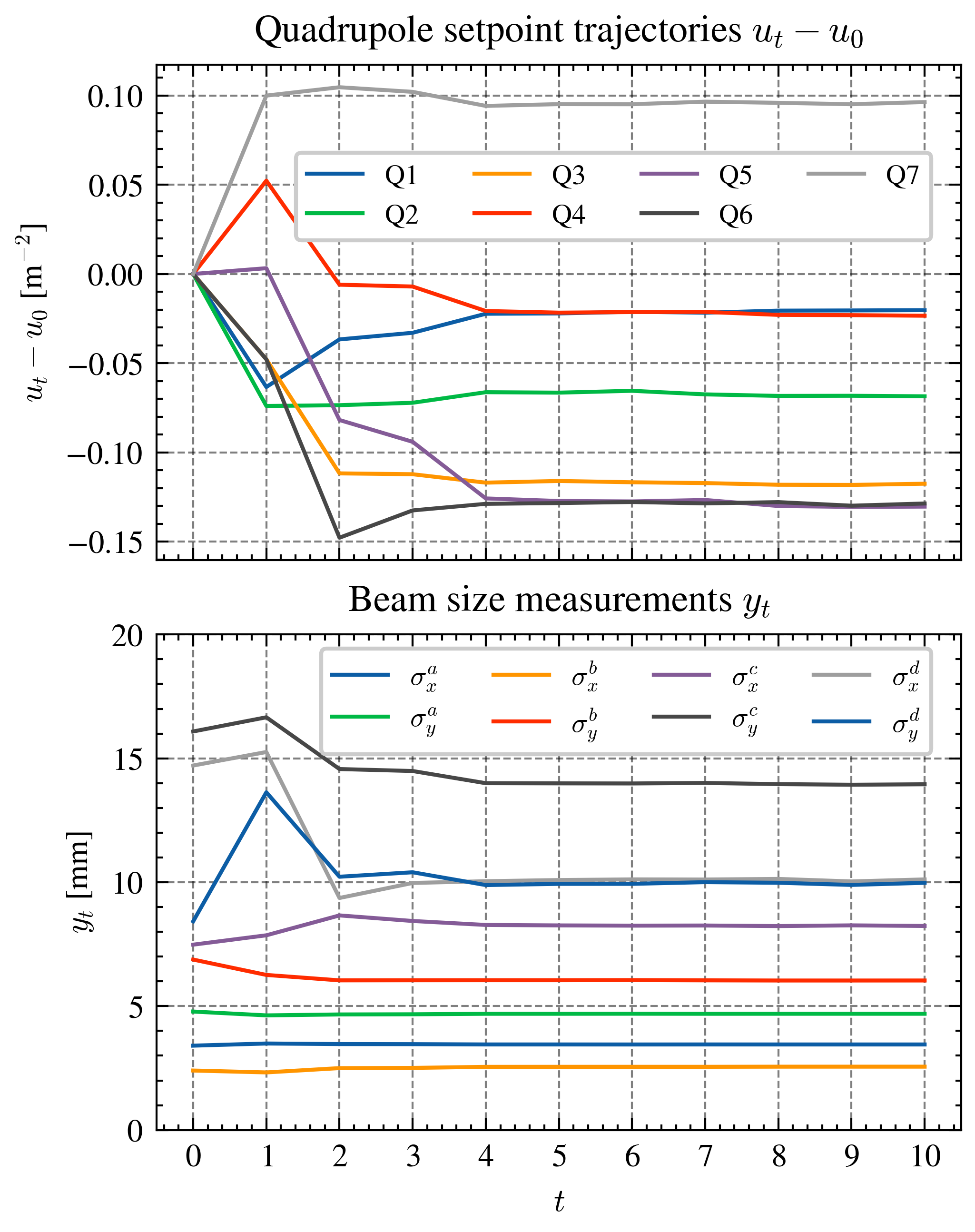}
\caption{Quadrupole setpoint values (top) and beam size measurements (bottom) over 11 time steps.}
\label{fig:results_accel_1}
\end{figure}

\begin{figure}[t]
\centering
\includegraphics[width=\columnwidth]{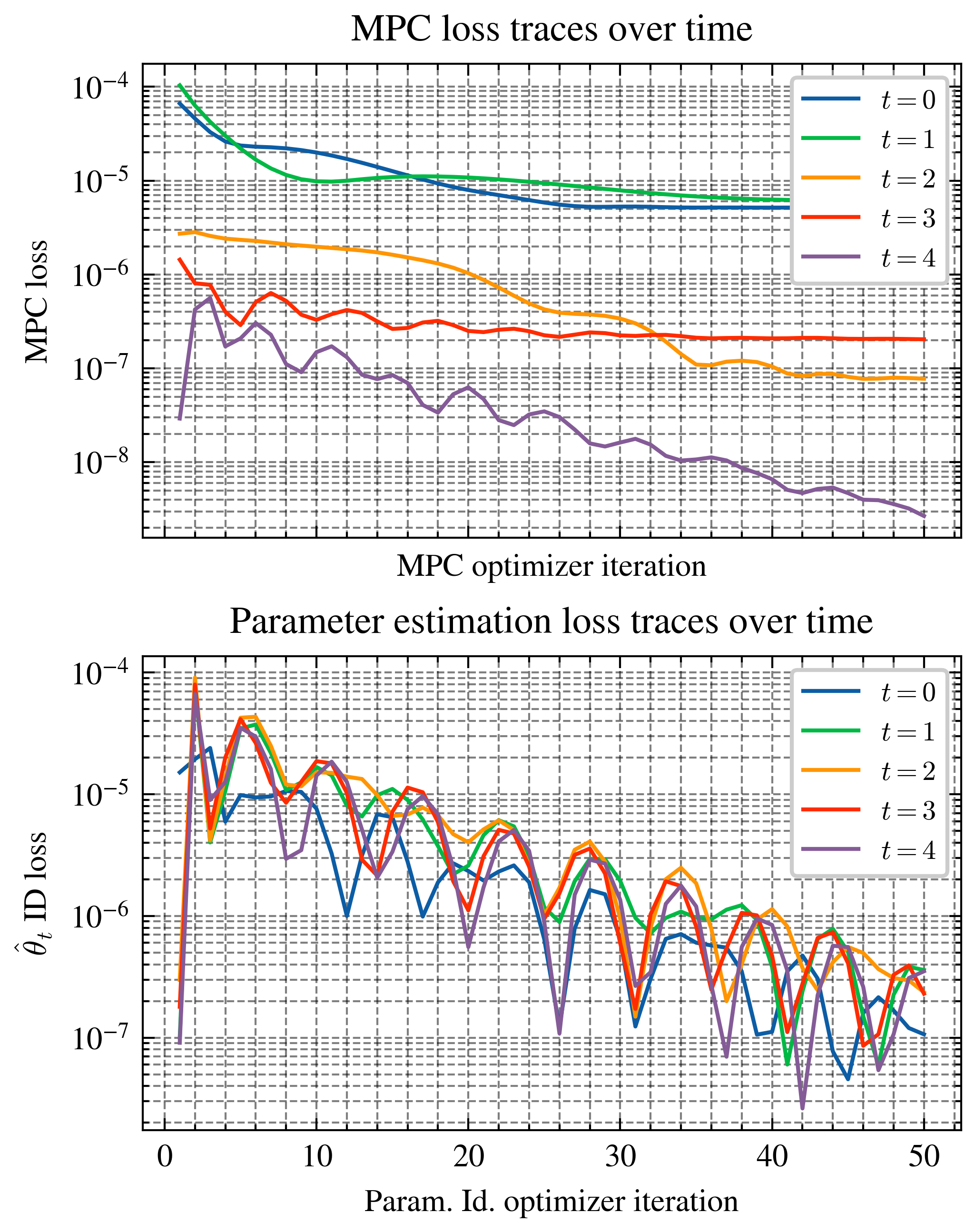}
\caption{Loss curves of \eqref{eq:ipf_cost_used} for the MPC optimization problem \eqref{eq:acc_static_mpc} (top) and \eqref{eq:acc_param_est} (bottom) for the parameter estimation for the first 5 time steps.}
\label{fig:results_accel_2}
\end{figure}

\subsection{Implementation}

At each time step $t$, we solve the static MPC problem \eqref{eq:acc_static_mpc} by differentiating through the Cheetah simulator and running Adam for 50 optimizer iterations with a learning rate of $10^{-2}$; at each iteration we roll out the horizon through the simulator and evaluate \eqref{eq:ipf_cost_used}. We warm-start the horizon sequence by shifting the previously optimized control sequence by one step (receding-horizon shift) and setting the new final input to the previous terminal value (i.e., $u_{t+N-1|t} \leftarrow u_{t+N-2|t-1}$). To enforce actuator feasibility \eqref{eq:acc_static_mpc_bounds}-\eqref{eq:acc_static_mpc_rate}, we evaluate the objective using a projected control sequence over the 50 optimizer steps: each planned input is clamped to a per-step change limit $\Delta u_{\max} = 0.1$ and to fixed sign/magnitude bounds relative to nominal setpoints. After optimization is complete, we apply only the first action $u_{t|t}$ to the plant and repeat.

Parameter estimation is performed using the least-squares objective \eqref{eq:acc_param_est} with $\lambda_\theta = 0$. In this example we use an expanding batch window: at control step $t$ we fit $\theta$ using all measurements collected so far, i.e., the dataset $\{(u_k,y_k)\}_{k=0}^{t}$, so the window size grows from 1 up to 11 steps over $t=0,\dots,10$. We update $\theta=[\sigma_x,\sigma_{p_x},\sigma_y,\sigma_{p_y}]^\top$ by running Adam for 50 iterations at each time step $t$ with a learning rate of $10^{-3}$. After each optimizer step, the parameters are projected onto the box constraint set $\Theta$. Instead of enforcing a hard constraint \eqref{eq:acc_static_mpc_g}, we include a hinge-squared soft constraint penalty
$\Omega_{\mathrm{soft}}$ in the objective.
The soft penalty maintains feasibility even when the initial beam size violates the diagnostic threshold $y_{\max}$,
as occurs in our simulation initial condition.

\subsection{Results}

Over 11 tuning iterations, the differentiable MPC controller drives the simulated IPF beam toward the desired circular spot at the target while respecting diagnostic safety limits. As shown in Fig.~\ref{fig:results_accel_1}, the run begins from an infeasible operating point in which two diagnostics violate the envelope threshold, with $\sigma_y^{c}$ and $\sigma_x^{d}$ initially exceeding $y_{\max}=14$~mm. The soft constraint penalty in \eqref{eq:ipf_soft_constraint} nevertheless keeps the optimization feasible and quickly steers the trajectory back into the safe region: within a few updates, all eight measurements fall below 14~mm and remain there, while the terminal diagnostic sizes $\sigma_x^{d}$ and $\sigma_y^{d}$ converge to the 10~mm target. The quadrupole setpoints evolve gradually (Fig.~\ref{fig:results_accel_1}, top), aided by the rate-limiting projection, which prevents overly aggressive corrective moves at the first update. Fig.~\ref{fig:results_accel_2} shows that the MPC problem \eqref{eq:acc_static_mpc} converges reliably within each control step (top), while the parameter estimation objective \eqref{eq:acc_param_est} rapidly decreases (bottom), producing a good fit after the first update. A small performance dip immediately after the first control action is expected because control is applied before system identification at $t=0$; once the first parameter update is incorporated, subsequent MPC steps improve both tracking at segment $d$ and satisfaction of the beam-size limits along the beamline.

\section{Conclusion}

We presented an MPC and online parameter estimation framework for systems whose dynamics are available only through a differentiable physics simulator. By differentiating through simulator rollouts, we optimize MPC objectives with respect to control inputs while updating unknown physical parameters from measurement mismatch using the same differentiable model. In the accelerator example, this enables gradient-based tuning of quadrupole setpoints while estimating an unknown incoming beam distribution; in the fluid example, the same approach supports control of a swimmer while estimating an unknown fluid parameter. Identifiability remains a key challenge when informative excitation is suppressed by closed-loop behavior. Future work will incorporate principled constrained solvers, excitation-aware objectives, and experiments utilizing state estimation for partially observed, high-dimensional systems.

\section*{Acknowledgments}
The authors thank Michael Brenner and Alex Scheinker for helpful discussions and guidance throughout this work.

\bibliographystyle{IEEEtranS}
\bibliography{references}
\end{document}